\documentclass[10pt]{article}
\usepackage{hyperref,amssymb,amsmath,graphicx,verbatim,amsfonts,amsthm}
\newtheorem{thm}{Theorem}[section]
\newtheorem{cor}[thm]{Corollary}
\newtheorem{lem}[thm]{Lemma}
\newtheorem{prop}[thm]{Proposition}
\newtheorem{clm}[thm]{Claim}

\newtheorem*{thm*}{Theorem}

\theoremstyle{definition}
\newtheorem{dfn}[thm]{Definition}
\theoremstyle{remark}

\numberwithin{equation}{section}

\newcommand{\set}[1]{\left\{#1\right\}}

\newcommand{\br}[1]{\left[#1\right]}
\newcommand{\Br}[2]{ \left[#1 \ \big| \ #2 \right] }

\newcommand{\Integer}{\mathbb{Z}}

\newcommand{\Z}{\Integer}
\newcommand{\N}{\mathbb{N}}

\newcommand{\eps}{\varepsilon}

\newcommand{\eqdef}{\stackrel{\mathrm{def}}{=}}

\DeclareMathOperator*{\E}{\mathbb{E}}

\renewcommand{\Pr}{}
\let\Pr\relax
\DeclareMathOperator*{\Pr}{\mathbb{P}}

\newcommand{\1}[1]{\mathbf{1}_{\set{ #1 } }}

\def\squareforqed{\hbox{\rlap{$\sqcap$}$\sqcup$}}
\def\qed{\ifmmode\squareforqed\else{\unskip\nobreak\hfil
\penalty50\hskip1em\null\nobreak\hfil\squareforqed
\parfillskip=0pt\finalhyphendemerits=0\endgraf}\fi}

\newcommand{\Cov}{\mathrm{Cov}}
\newcommand{\g}{\gamma}

\renewcommand{\1}[1]{\mathbf{1}{\set{ #1 } }}

\renewcommand{\o}{\omega}

\newcommand{\Ups}{\Upsilon}
\newcommand{\ups}{\upsilon}

\newcommand{\vphi}{\varphi}

\newcommand{\id}{\mathbf{id}}

\newcommand{\supp}{\mathrm{supp}}

\newcommand{\Ee}{\mathcal{E}}
\begin{document}

\title{Rate of Escape of the Mixer Chain}

\author{Ariel Yadin\thanks{Faculty of Mathematics and Computer Science, The
Weizmann Institute of Science, Rehovot 76100, Israel.
Email: \texttt{ariel.yadin@weizmann.ac.il}} }

\date{ }

\maketitle

\begin{abstract}
The mixer chain on a graph $G$ is the following Markov chain.
Place tiles on the vertices of $G$, each tile labeled by its
corresponding vertex. A ``mixer'' moves randomly on the graph,
at each step either moving to a randomly chosen neighbor, or swapping
the tile at its current position with some randomly chosen adjacent tile.

We study the mixer chain on $\Z$, and show that at time $t$
the expected distance to the origin is $t^{3/4}$, up to constants.
This is a new example of a random walk
on a group with rate of escape strictly between $t^{1/2}$ and $t$.
\end{abstract}

\section{Introduction}

Let $G = (V,E)$ be a graph.  On each vertex $v \in V$, place a tile marked $v$.
Consider the following Markov chain, which we call the \emph{mixer chain}.  A ``mixer''
performs a random walk on the graph.  At each time step, the mixer chooses a random vertex
adjacent to its current position.  Then, with probability $1/2$ it moves to that vertex, and with
probability $1/2$ it remains at the current location, but swaps the tiles on the current vertex and the adjacent
vertex.  If $G$ is the Cayley graph of a group, then the mixer chain turns out to be
a random walk on a different group.

Aside from being a canonical process, the mixer chain is interesting because of its \emph{rate of escape}.
For a random walk $\set{X_t}$ on some graph $G$, we use the terminology \emph{rate of escape} for the limit
$$ \lim_{t \to \infty} \frac{\log \E [d(X_t,X_0)]}{\log t} , $$
where $d(\cdot,\cdot)$ is the graphical distance.
When restricting to random walks on groups, it is still open what values in $[0,1]$ can be obtained
by rates of escape.  For example, if the group is $\Z^d$ then the rate of escape is $1/2$.
On a $d$-ary tree (free group) the rate of escape is $1$.
As far as the author is aware,
the only other examples known were given by Erschler in \cite{Erschler} (see also \cite{Revelle}).
Erschler iterates a construction known as the lamp-lighter (slightly similar to the mixer chain),
and produces examples of groups with rates of escape $1-2^{-k}$, $k=1,2,\ldots,$.

After formally defining the mixer chain on general groups, we study the mixer chain on $\Z$.
Our main result, Theorem \ref{thm:main thm}, shows that the mixer chain on $\Z$ has rate
of escape $3/4$.

It is not difficult to show (perhaps using ideas from this note)
that on transient groups the mixer chain has rate of escape $1$.
Since all recurrent groups are essentially $\Z$ and $\Z^2$,
it seems that the mixer chain on other groups cannot give examples of other rates of escape.
As for $\Z^2$, one can show that the mixer chain has rate of escape $1$.
In fact, the ideas in this note suggest that the distance to the origin in the mixer chain on $\Z^2$ is $n \log^{-1/2}(n)$ up to constants.

After introducing some notation,
we provide a formal definition of the mixer chain, as random walk on
a Cayley graph.  The generalization to general graphs is immediate.

\paragraph{Acknowledgement.}
I wish to thank Itai Benjamini for suggesting this construction,
and for useful discussions.

\subsection{Notation}
\label{scn:notation}

Let $G$ be a group and $U$ a generating set for $G$, such that if
$x \in U$ then $x^{-1} \in U$ ($U$ is called \emph{symmetric}).
The \emph{Cayley graph} of $G$ with respect to $U$
is the graph with vertex set $G$ and edge set
$\set{ \set{g,h} \ : \  g^{-1} h \in U }$.
Let ${\cal D}$ be a distribution on $U$.  Then we can define
the \emph{random walk} on $G$ (with respect to $U$ and ${\cal D}$)
as the Markov chain with state space $G$ and transition matrix
$P(g,h) = \1{ g^{-1} h \in U } {\cal D}(g^{-1}h)$.
We follow the convention that such a process starts from the
identity element in $G$.

A permutation of $G$ is a bijection from $G$ to $G$.
The \emph{support} of a permutation $\sigma$, denoted $\supp(\sigma)$,
is the set of all elements $g \in G$
such that $\sigma(g) \neq g$.
Let $\Sigma$ be the group of all permutations of $G$ with finite support
(multiplication is composition of functions).
By $<g,h>$ we denote the transposition of $g$ and $h$; that is, the permutation
$\sigma$ with support $\set{g,h}$ such that $\sigma(g) = h$, $\sigma(h) = g$.
By $<g_1,g_2,\ldots,g_n>$ we denote the cyclic permutation
$\sigma$ with support $\set{g_1,\ldots,g_n}$, such that $\sigma(g_j) = g_{j+1}$
for $j<n$ and $\sigma(g_n) = g_1$.

For an element $g \in G$ we associate a canonical permutation,
denoted by $\phi_g$, defined by $\phi_g(h) = gh$ for all $h \in G$.
It is straightforward to verify that the map $g \mapsto \phi_g$ is a homomorphism
of groups, and so we use $g$ to denote $\phi_g$.
Although $g \not\in \Sigma$, we have that $g \sigma g^{-1} \in \Sigma$ for all
$\sigma \in \Sigma$.

We now define a new group, that is in fact the
\emph{semi-direct product} of $G$ and $\Sigma$, with respect to
the homomorphism $g \mapsto \phi_g$ mentioned above.
The group is denoted by $G \ltimes \Sigma$, and its elements are $G \times \Sigma$.
Group multiplication is defined by:
$$ (g,\sigma) (h, \tau) \eqdef (gh, g \tau g^{-1} \sigma  ) . $$
We leave it to the reader to verify that this is a well-defined
group operation.
Note that the identity element in this group is $(e, \id)$, where $\id$
is the identity permutation in $\Sigma$ and $e$ is the identity element in $G$.
Also, the inverse of $(g,\sigma)$ is $(g^{-1}, g^{-1} \sigma^{-1} g)$.

%
%
%

We use $d(g,h) = d_{G,U}(g,h)$ to denote the distance between $g$ and $h$ in the
group $G$ with respect to the generating set $U$; i.e., the minimal $k$ such that
$g^{-1}h = \prod_{j=1}^k u_j$ for some $u_1,\ldots,u_k \in U$.
The generating set also provides us with a graph structure.
$g$ and $h$ are said to be adjacent if $d(g,h) = 1$,
that is if $g^{-1} h \in U$.
A path $\g$ in $G$ (with respect to the generating set $U$) is a sequence
$(\g_0,\g_1,\ldots,\g_n)$.  $|\g|$ denotes the \emph{length} of the path,
which is defined as the length of the sequence minus $1$ (in this case $|\g| = n$).

\subsection{Mixer Chain}

In order to define the mixer chain we require the following
\begin{prop} \label{prop:generating set}
Let $U$ be a finite symmetric generating set for $G$.  Then,
$$ \Upsilon = \set{ (u,\id) , (e,<e,u>) \ : \ u \in U } $$
generates $G \ltimes \Sigma$.  Furthermore,
for any cyclic permutation $\sigma = <g_1,\ldots,g_n> \in \Sigma$,
$$ d_{G \ltimes \Sigma, \Upsilon} ( (g_1,\sigma) , (g_1,\id) )
\leq 5 \sum_{j=1}^{n} d(g_j,\sigma(g_j)) . $$
\end{prop}

\begin{proof}
Let $D((g,\sigma),(h,\tau))$ denote the minimal $k$ such that
$(g,\sigma)^{-1} (h,\tau) = \prod_{j=1}^k \ups_j$, for some $\ups_1,\ldots,\ups_k \in \Ups$,
with the convention that $D((g,\sigma), (h,\tau)) = \infty$ if there is no such finite sequence of
elements of $\Ups$.  Thus, we want to prove that $D((g,\sigma),(e,\id)) < \infty$ for all $g \in G$
and $\sigma \in \Sigma$.
Note that by definition
for any $f \in G$ and $\pi \in \Sigma$,
$D((g,\sigma),(h,\tau)) = D((f,\pi)(g,\sigma),(f,\pi)(h,\tau))$.

A \emph{generator simple path} in $G$ is a finite sequence of generators
$u_1,\ldots,u_k \in U$ such that for any
$1 \leq \ell \leq k$, $\prod_{j=\ell}^k u_j \neq e$.
By induction on $k$, one can show that for any $k \geq 1$,
and for any generator simple path $u_1,\ldots, u_k$,
\begin{align} \label{eqn:generating product}
(e, <e, \prod_{j=1}^k u_j > ) = \prod_{j=1}^{k-1} (e, <e,u_j>) (u_j,\id) \cdot
(e, <e,u_k>) \cdot \prod_{j=1}^{k-1} (e, <e,u_{k-j}^{-1}>) (u_{k-j}^{-1},\id) .
\end{align}
If $d(g,h) = k$ then there exists a generator simple path
$u_1,\ldots,u_k$ such that $h = g \prod_{j=1}^k u_j$.
Thus, we get that for any $h \in G$,
$$ D((e,<e,h>),(e,\id)) \leq 4 d(h,e) - 3 . $$
Because $g<e,g^{-1}h>g^{-1} = <g,h>$, we get that if $\tau = <g,h>\sigma$ then
$$ D((g,\tau),(g,\sigma)) = D((g,\sigma)(e,<e,g^{-1}h>),(g,\sigma)(e,\id))
\leq 4 d(g^{-1}h,e) - 3 = 4 d(g,h) - 3 . $$
The triangle inequality now implies that
$D((h,\tau),(g,\sigma)) \leq 5 d(g,h) - 3$.

Thus, if $\sigma = <g_1,g_2,\ldots,g_n>$, since $\sigma = <g_1,g_2><g_2,g_3> \cdots <g_{n-1},g_n>$,
we get that
\begin{align} \label{eqn:D(cycle)}
D((g_1,\sigma), (g_1,\id)) \leq 5 \sum_{j=1}^{n-1} d(g_j,g_{j+1}) + d(g_n,g_1) .
\end{align}
The proposition now follows from the fact that any $\sigma \in \Sigma$ can be written as a finite
product of cyclic permutations.
\end{proof}

We are now ready to define the mixer chain:
\begin{dfn} \label{dfn:mixer}
Let $G$ be a group with finite symmetric generating set $U$.
The mixer chain on $G$ (with respect to $U$)
is the random walk on the group $G \ltimes \Sigma$ with respect to uniform measure
on the generating set $\Upsilon = \set{ (u,\id) , (e,<e,u>) \ : \ u \in U }$.
\end{dfn}

An equivalent way of viewing this chain is viewing the state $(g,
\sigma) \in G \ltimes \Sigma$ as follows:  The first coordinate
corresponds to the position of the mixer on $G$.
The second coordinate corresponds to the placing of the different
tiles, so the tile marked $x$ is placed on the vertex $\sigma(x)$.
By Definition \ref{dfn:mixer},
the mixer chooses uniformly an adjacent vertex of $G$, say $h$.  Then,
with probability $1/2$ the mixer swaps the tiles on $h$ and $g$,
and with probability $1/2$ it moves to $h$.  The identity element in
$G \ltimes \Sigma$ is $(e,\id)$, so the mixer starts at $e$ with all tiles
on their corresponding vertices (the identity permutation).

\subsection{Distance Bounds}

In this section we show that the distance of an element in $G \ltimes \Sigma$
to $(e,\id)$ is essentially governed by the sum of
the distances of each individual tile to its origin.

Let $(g,\sigma) \in G \ltimes \Sigma$.
Let $\gamma = (\g_0,\g_1,\ldots,\g_n)$ be a finite path
in $G$. 
We say that the path $\gamma$ \emph{covers} $\sigma$ if
$\supp(\sigma) \subset \set{\g_0,\g_1,\ldots,\g_n}$.
The \emph{covering number} of $g$ and $\sigma$, denoted $\Cov(g,\sigma)$,
is the minimal length of a path $\gamma$, starting at $g$, that covers $\sigma$;
i.e.
$$ \Cov(g,\sigma) = \min \set{ |\g| \ : \
\gamma_0 = g \textrm{ and } \gamma \textrm{ is a path covering } \sigma } . $$

To simplify notation, we denote $D = d_{G \ltimes \Sigma, \Ups}$.

\begin{prop} \label{prop:distance upper in mixer group}
Let $(g,\sigma) \in G \ltimes \Sigma$.
Then,
$$ D((g,\sigma),(g,\id)) \leq
2 \Cov(g,\sigma) + 5 \sum_{h \in \supp(\sigma)} d(h,\sigma(h)) . $$
\end{prop}

\begin{proof}
The proof of the proposition is by induction on the size of $\supp(\sigma)$.
If $|\supp(\sigma)| = 0$, then $\sigma = \id$
so the proposition holds.
Assume that $|\supp(\sigma)| > 0$.

Let $n  = \Cov(g,\sigma)$,
and let $\g$ be a path in $G$ such that $|\g| = n$,
$\g_0 = g$ and $\g$ covers $\sigma$.
Write $\sigma = c_1 c_2 \cdots c_k$, where the $c_j$'s
are cyclic permutations with pairwise disjoint non-empty supports, and
$$ \supp(\sigma) = \bigcup_{j=1}^k \supp(c_j) . $$

Let
$$ j = \min \set{ m \geq 0 \ : \ \g_m \in \supp(\sigma) } . $$
So, there is a unique $1 \leq \ell \leq k$ such that
$\g_j \in \supp(c_\ell)$.
Let $\tau = c_{\ell}^{-1} \sigma$.  Thus,
$$ \supp(\tau) = \bigcup_{j \neq \ell} \supp(c_j) , $$
and specifically, $|\supp(\tau)| < |\supp(\sigma)|$.
Note that $h \in \supp(\g_j^{-1}c_\ell \g_j)$ if and only if
$\g_j h \in \supp(c_\ell)$,
and specifically, $e \in \supp(\g_j^{-1} c_\ell \g_j)$.
$\g_j^{-1} c_\ell \g_j$ is a cyclic permutation, so
by Proposition \ref{prop:generating set}, we know that
\begin{align} \label{eqn:induction step D(cycles)}
D((\g_j,\sigma),(\g_j,\tau)) & = D((\g_j,\tau)(e,\g_j^{-1} c_\ell \g_j), (\g_j,\tau))
= D((e, \g_j^{-1} c_\ell \g_j), (e,\id)) \nonumber \\
& \leq 5 \sum_{h \in \supp(c_\ell)} d(\g_j^{-1} h, \g_j^{-1} c_\ell (h))
= 5 \sum_{h \in \supp(c_\ell)} d(h,\sigma(h)) .
\end{align}
By induction,
\begin{align} \label{eqn:tau in induction step}
D((\g_j,\tau),(\g_j,\id)) \leq
2 \Cov(\gamma_j,\tau)
+ 5 \sum_{h \in \supp(\tau)} d(h,\tau(h)) .
\end{align}
Let $\beta$ be the path $(\g_j,\g_{j+1}, \ldots, \g_n)$.
Since $\g_j$ is the first element in $\g$ that is in $\supp(\sigma)$,
we get that $\supp(\tau) \subset \supp(\sigma) \subseteq
\set{ \g_j,\g_{j+1}, \ldots, \g_n }$,
which implies that $\beta$ is a path of length $n-j$ that covers $\tau$,
so $\Cov(\g_j,\tau) \leq n-j$.
Combining \eqref{eqn:induction step D(cycles)} and \eqref{eqn:tau in induction step} we get,
\begin{align*}
D((g,\sigma),(g,\id)) & \leq D((\g_0,\sigma),(\g_j,\sigma)) +
D((\g_j,\sigma),(\g_j,\tau)) + D((\g_j,\tau),(\g_j,\id))
+ D((\g_j,\id),(\g_0,\id)) \\
& \leq j + 5 \sum_{h \in \supp(c_\ell)} d(h,\sigma(h))
+ 5 \sum_{h \in \supp(\tau)} d(h,\sigma(h))
+ 2 (n-j) + j \\
& = 2 \Cov(g,\sigma) + 5 \sum_{h \in \supp(\sigma)} d(h,\sigma(h)) .
\end{align*}
\end{proof}

\begin{prop}\label{prop:distance lower in mixer group}
Let $(g,\sigma) \in G \ltimes \Sigma$ and let $g' \in G$.  Then,
$$ D((g,\sigma),(g',\id)) \geq \frac{1}{2} \sum_{h \in \supp(\sigma)} d(h,\sigma(h)) . $$
\end{prop}

\begin{proof}
The proof is by induction on $D = D((g,\sigma),(g',\id))$.
If $D = 0$ then $\sigma = \id$, and we are done.
Assume that $D>0$.  Let $\ups \in \Ups$ be a generator
such that $D((g,\sigma) \ups ,(g',\id)) = D-1$.
There exists $u \in U$ such that either
$\ups = (u,\id)$ or $\ups = (e,<e,u>)$.
If $\ups = (u,\id)$ then by induction
$$ D \geq D((g,\sigma) \ups, (g',\id)) \geq \frac{1}{2} \sum_{h \in \supp(\sigma)} d(h,\sigma(h)) .  $$
So assume that $\ups = (e,<e,u>)$.
If $\sigma(h) \not\in \set{g,gu}$, then $<g,gu> \sigma (h) = \sigma (h)$,
and
$$ \supp(\sigma) \setminus \set{ \sigma^{-1}(g), \sigma^{-1}(gu) } = \supp(<g,gu> \sigma) \setminus
\set{ \sigma^{-1}(g), \sigma^{-1}(gu) } . $$
Since $d(g,gu) = 1$,
\begin{align*}
\sum_{h \in \supp(\sigma)} d(h,\sigma(h)) & =
d(g,\sigma^{-1}(g)) + d(gu,\sigma^{-1}(gu)) + \sum_{h \not\in \set{\sigma^{-1}(g),\sigma^{-1}(gu)}} d(h,\sigma(h)) \\
& \leq d(g,gu) + d(gu, \sigma^{-1}(g)) + d(gu,g) + d(g,\sigma^{-1}(gu)) \\
& \qquad +
\sum_{h \not\in \set{\sigma^{-1}(g),\sigma^{-1}(gu)}} d(h,<g,gu>\sigma(h))  \\
& \leq 2 + \sum_{h \in \supp(<g,gu>\sigma)} d(h,<g,gu>\sigma(h)) .
\end{align*}
So by induction,
\begin{align*}
D & = 1 + D((g,<g,gu>\sigma) , (g',\id)) \geq 1 + \frac{1}{2} \sum_{h \in \supp(<g,gu>\sigma) } d(h,<g,gu>\sigma(h)) \\
& \geq \frac{1}{2} \sum_{h \in \supp(\sigma)} d(h,\sigma(h))  .
\end{align*}
\end{proof}

\section{The Mixer Chain on $\Z$}

We now consider the mixer chain on $\Z$, with $\set{1,-1}$ as the symmetric generating set.
We denote by $\set{\o_t = (S_t,\sigma_t)}_{t \geq 0}$ the mixer chain on $\Z$.

For $\o \in \Z \ltimes \Sigma$ we denote by $D(\o)$ the distance
of $\o$ from $(0,\id)$ (with respect to the generating set $\Ups$, see
Definition \ref{dfn:mixer}).
Denote by $D_t = D(\o_t)$ the distance of the
chain at time $t$ from the origin.

As stated above, we show that the mixer chain on $\Z$ has rate of escape $3/4$.
In fact, we prove slightly stronger bounds on the distance to the origin at time $t$.

\begin{thm} \label{thm:main thm}
Let $D_t$ be the distance to the origin of the mixer chain on $\Z$.
Then, there exist constants $c,C>0$ such that for all $t \geq 0$,
$c t^{3/4} \leq \E[D_t] \leq C t^{3/4}$.
\end{thm}

The proof of Theorem \ref{thm:main thm} is in Section \ref{sec:proof of thm}.

For $z \in \Z$, denote by $X_t(z) = |\sigma_t(z)-z|$,
the distance of the tile marked $z$
to its origin at time $t$.
Define
$$ X_t = \sum_{z \in \Z} X_t(z) , $$
which is a finite sum for any given $t$.
As shown in Propositions \ref{prop:distance upper in mixer group} and
\ref{prop:distance lower in mixer group},
$X_t$ approximates $D_t$ up to certain factors.

For $z \in \Z$ define
$$ V_t(z) = \sum_{j=0}^t \1{ S_t = \sigma_t(z) } . $$
$V_t(z)$ is the number of times
that the mixer visits the tile marked $z$,
up to time $t$.

\subsection{Distribution of $X_t(z)$}

\begin{prop} \label{prop:(z) and (-z)}
For $\sigma \in \Sigma$ define $\sigma' \in \Sigma$ by
$\sigma'(z) = -\sigma(-z)$ for all $z \in \Z$.
Then, for any $t \geq 1$,
$((S_1,\sigma_1),\ldots,(S_t,\sigma_t))$ and $((-S_1,\sigma'_1),\ldots,(-S_t,\sigma'_t))$
have the same distribution.
\end{prop}

\begin{proof}
Let $\vphi$ be the permutation (not in $\Sigma$) defined by $\vphi(x) = -x$ for all $x \in \Z$.
Then, $\sigma' = \vphi \sigma \vphi$.
Since $\vphi^2 = \id$, we get that $(\sigma \tau)' = \sigma' \tau'$, $\sigma'' = \sigma$ and
$(\sigma^{-1})' = (\sigma')^{-1}$.

The proof is by induction on $t$.
If $t = 1$, then $(S_1,\sigma_1)$
is uniformly distributed over the set
$$ \Ups = \set{ (1,\id) , (-1,\id) , (0,<0,1>) , (0,<0,-1>) } . $$
Since $<0,1>' = <0,-1>$, and $\id' = \id$,
the proposition is proved for $t=1$.

Let $t>1$.
Let $(x,\tau) \in \Z \ltimes \Sigma$ be any element,
and let $\ups = (y,\rho) \in \Ups$ be a generator.
We have
$(x,\tau)(y,\rho) = (x+y, x \rho (-x) \tau)$.
Since $\rho \in \set{ \id, <0,1>, <0,-1>}$,
we have that
$(x \rho (-x) \tau)' = (-x) \rho' x \tau'$.
Since $y \neq 0$ if and only if $\rho = \id$,
we get that $(-y) \rho y = \rho$, so
$$ (-x,\tau')(-y,\rho') =
(-(x+y), (-x) \rho' x \tau') =
(-(x+y), (x \rho (-x) \tau)' ) . $$
Thus, we get that for any $(z,\sigma),(x,\tau) \in \Z \ltimes \Sigma$ and any
$(y,\rho) \in \Ups$,
$$ (z,\sigma) = (x,\tau) (y,\rho) \qquad \textrm {if and only if } \qquad
(-z,\sigma') = (-x,\tau') (-y,\rho') . $$
If $(y,\rho)$ is uniformly distributed in $\Ups$, then so is $(-y,\rho')$.
Thus, since $(S_{t-1},\sigma_{t-1})^{-1} (S_t,\sigma_t)$ is uniformly distributed in $\Ups$,
$(-S_{t-1},\sigma_{t-1}')^{-1} (-S_t,\sigma_{t}')$ is also uniformly distributed in $\Ups$.
By induction,
$((S_1,\sigma_1),\ldots,(S_{t-1},\sigma_{t-1}))$ and
$((-S_1,\sigma'_1),\ldots,(-S_{t-1},\sigma'_{t-1}))$, have the same distribution.
The Markov property now implies the proposition.
\end{proof}

By a lazy random walk on $\Z$, we refer to the integer valued process $W_t$,
such that $W_{t+1} - W_t$ are i.i.d. random variables with the distribution
$\Pr [ W_{t+1} - W_t = 1 ] = \Pr [ W_{t+1} - W_t = 1 ] = 1/4$ and
$\Pr [ W_{t+1} - W_t = 0 ] =1/2$.

\begin{lem} \label{lem:conditioned on number of visits}
Let $t \geq 0$ and $z \in \Z$.  Let $k \geq 1$ such that $\Pr [V_t(z)=k] > 0$.
Then,
conditioned on $V_t(z) = k$, the distribution of $\sigma_t(z)-z$ is
the same as $W_{k-1} + B$,
where $\set{W_k}$ is a lazy random walk on $\Z$, and $B$ is a random variable
independent of $\set{W_k}$ such that $|B| \leq 2$.
\end{lem}

\begin{proof}
Define inductively the following random times:
$T_0(z) = 0$, and for $j \geq 1$,
$$ T_j(z) = \inf \set{ t \geq T_{j-1}(z)+1 \ : \ S_t = \sigma_t(z) } . $$

\begin{clm} \label{clm:sigma_T}
Let $T = T_1(0)$.
For all $\ell$ such that $\Pr [T=\ell] > 0$,
$$ \Pr \Br{ \sigma_T(0) =  1 }{ T = \ell } =
\Pr \Br{ \sigma_T(0) = - 1 }{ T = \ell } = 1/4 , $$
and
$$ \Pr \Br{ \sigma_T(0) = 0 }{ T = \ell } = 1/2 , $$
\end{clm}

\begin{proof}
Note that $|S_1-\sigma_1(0)| = 1$ and that for all $1 \leq t < T$,
$\sigma_t(0) = \sigma_1(0)$.  Thus, $\sigma_{T-1}(0) = \sigma_1(0)$ and $S_{T-1} = S_1$.
So we have the equality of events
$$ \set{ T=\ell} = \bigcap_{t=1}^{\ell-1} \set{ S_t \neq \sigma_t(0) } \bigcap
\set{ S_{\ell-1}=S_1 , \sigma_{\ell-1}(0) = \sigma_1(0) }
\bigcap \set{ S_\ell = \sigma_1(0) \textrm{ or } \sigma_\ell(0) = S_1 } . $$
Hence, if we denote $\Ee = \bigcap_{t=1}^{\ell-1} \set{ S_t \neq \sigma_t(0) } \bigcap
\set{ S_{\ell-1}=S_1 , \sigma_{\ell-1}(0) = \sigma_1(0) }$, then
\begin{align} \label{eqn:prob of T=ell}
\Pr [T=\ell] & = \Pr [ \Ee ] \cdot
\Pr \Br{ S_\ell = \sigma_1(0) \textrm{ or } \sigma_\ell(0) = S_1 }{ \Ee } \nonumber \\
& = \Pr [ \Ee] \cdot \frac{1}{2} .
\end{align}
Since the events $\set{S_1=0}$ and $\set{\sigma_1(0)=0}$ are disjoint and their union
is the whole space, we get that
\begin{align} \label{eqn:prob sigma=0,T=ell}
\Pr [ \sigma_T(0) = 0 , T=\ell ] & = \Pr [ \Ee , S_\ell = \sigma_1(0) = 0 ] +
\Pr [ \Ee , \sigma_\ell(0) = S_1 = 0 ] \nonumber \\
& = \Pr \br{ \Ee , \sigma_1(0) = 0 } \cdot
\Pr \Br{ S_{\ell} = \sigma_1(0) }{ S_{\ell-1} = S_1 , \sigma_{\ell-1}(0) = \sigma_1(0) = 0 } \nonumber \\
& \quad + \Pr \br{ \Ee , S_1(0) = 0 } \cdot
\Pr \Br{ \sigma_{\ell}(0) = S_1 }{ S_{\ell-1} = S_1 = 0 , \sigma_{\ell-1}(0) = \sigma_1(0) } \nonumber  \\
& = \Pr [\Ee] \cdot \frac{1}{4} .
\end{align}
Combining \eqref{eqn:prob of T=ell} and \eqref{eqn:prob sigma=0,T=ell} we get that
$$ \Pr \Br{ \sigma_T(0)=0}{ T=\ell } = \frac{1}{2} . $$
Finally, by Proposition \ref{prop:(z) and (-z)},
\begin{align*}
\Pr \br{ \sigma_T(0) = 1 , T = \ell } & = \Pr \br{ \Ee \ , \
S_\ell = \sigma_\ell(0) = 1 } \\
& = \Pr \br{ \sigma_T(0) = -1, T=\ell } .
\end{align*}
Since the possible values for $\sigma_T(0)$ are $-1,0,1$, the claim follows.
\end{proof}

%
%
%

%
%

We continue with the proof of Lemma \ref{lem:conditioned on number of visits}.

We have the equality of events $\set{V_t(z) = k } =
\set{ T_k(z) \leq t < T_{k+1}(z) }$.

Let $t_1,t_2,\ldots,t_k,t_{k+1}$ be such that
$$ \Pr [ T_1(z)=t_1, \ldots, T_{k+1}(z) = t_{k+1} ] > 0 , $$
and condition on the event
$\Ee = \set{ T_1(z)=t_1, \ldots, T_{k+1}(z) = t_{k+1} }$.
Assume further that $t_k \leq t < t_{k+1}$, so that $V_t(z)=k$.
Write
\begin{align} \label{eqn:sigma(z)-z = sum}
\sigma_t(z)-z = \sigma_t(z) - \sigma_{T_k(z)}(z) + \sum_{j=2}^{k}
\sigma_{T_j(z)}(z)-\sigma_{T_{j-1}(z)}(z) + \sigma_{T_1(z)}(z) - z .
\end{align}

For $1 \leq j \leq k-1$ denote $Y_j = \sigma_{T_{j+1}(z)}(z) - \sigma_{T_j(z)}(z)$.
By Claim \ref{clm:sigma_T} and the Markov property,
conditioned on $\Ee$, $\set{Y_j}$ are independent with the distribution
$\Pr[Y_j=1|\Ee] = \Pr[Y_j=-1|\Ee] = 1/4$ and $\Pr[Y_j=0|\Ee] = 1/2$.
So conditioned on $\Ee$, $\sum_{j=1}^{k-1} Y_j$ has the same distribution of
$W_{k-1}$.

Finally, $|\sigma_t(z) - \sigma_{T_k(z)}(z)| \leq 1$,
and $|\sigma_{T_1(z)}(z) - z| \leq 1$.
Since conditioned on $\Ee$, $\sigma_t(z) - \sigma_{T_k(z)}(z)$,
and $\sigma_{T_1(z)}(z) - z$ are independent of $\set{Y_j}$, this
completes the proof of the lemma.
\end{proof}

\begin{cor} \label{cor:E[Xt(z)]}
There exist constants $c,C>0$ such that for all $t \geq 0$
and all $z \in \Z$,
$$ c \E[ \sqrt{V_t(z)} ] - 2 \Pr [V_t(z) \geq 1] \leq \E [X_t(z)] \leq C \E[ \sqrt{V_t(z)} ] + 2 \Pr [V_t(z) \geq 1] . $$
\end{cor}

\begin{proof}
Let $\set{W_t}$ be a lazy random walk on $\Z$.  Note that
$\set{2 W_t}$ has the same distribution as $\set{S'_{2t}}$
where $\set{S'_t}$ is a simple random walk on $\Z$.
It is well known (see e.g. \cite{Revesz}),
that there exist universal constants $c_1,C_1>0$ such that
for all $t \geq 0$,
$$ c_1 \sqrt{t} \leq \E[|S'_{2t}|] = 2 \E [|W_t|] \leq C_1 \sqrt{t} . $$
By Lemma \ref{lem:conditioned on number of visits}, we know that
for any $k \geq 0$,
$$ \E[ |W_k| ] - 2 \leq \E \Br{ X_t(z) }{ V_t(z) =k+1 } \leq \E[ |W_k|] + 2 . $$
Thus, summing over all $k$, there exists constants $c_2,C_2>0$ such that
$$ c_2 \E[ \sqrt{V_t(z)} ] - 2\Pr [V_t(z) \geq 1] \leq \E [X_t(z)]
\leq C_2 \E[ \sqrt{V_t(z)} ] + 2\Pr [V_t(z) \geq 1] . $$
\end{proof}

\begin{lem} \label{lem:distribution of Vt(z)}
Let $\set{S'_t}$ be a simple random walk on $\Z$ started at $S'_0=0$,
and let
$$ L_t(z) = \sum_{j=0}^t \1{ S'_j=z } . $$
Then, for any $z \in \Z$, and any $k \in \N$,
$$ \Pr [ L_{2t}(2z) \geq k ] \leq \Pr [ V_t(z) \geq k ] . $$
Specifically, $\E[\sqrt{L_{2t}(2z)}] \leq \E[\sqrt{V_t(z)}]$.
\end{lem}

\begin{proof}
Fix $z \in \Z$.
For $t \geq 0$ define $M_t = S_t-\sigma_t(z)+z$.
Note that
$$ V_t(z) = \sum_{j=0}^t \1{ M_j = z} , $$
so $V_t(z)$ is the number of times $\set{M_t}$ visits $z$ up to time $t$.

$\set{M_t}$ is a Markov chain on $\Z$ with the following step distribution.
$$ \Pr \Br{ M_{t+1} = M_t + \eps }{ M_t } =
\left\{ \begin{array}{lrl}
1/2 & M_t = z ,& \eps \in \set{-1,1} , \\
1/2 & |M_t-z|=1 ,& \eps = - M_t+z , \\
1/4 & |M_t-z|=1 ,& \eps = 0 , \\
1/4 & |M_t-z|=1 ,& \eps = M_t-z , \\
1/4 & |M_t-z| > 1 ,& \eps \in \set{-1,1} , \\
1/2 & |M_t-z| > 1 ,& \eps = 0 .
\end{array} \right. $$
Specifically, $\set{M_t}$ is simple symmetric when at $z$, lazy symmetric when not adjacent to $z$,
and has a drift towards $z$ when adjacent to $z$.

Define $\set{N_t}$ to be the following Markov chain on $\Z$:
$N_0=0$, and for all $t \geq 0$,
$$ \Pr \Br{ N_{t+1} = N_t + \eps }{ N_t } =
\left\{ \begin{array}{lrl}
1/2 & N_t = z ,& \eps \in \set{-1,1} , \\
1/2 & N_t \neq z , & \eps = 0 , \\
1/4 & N_t \neq z , & \eps \in \set{-1,1} . \\
\end{array} \right. $$
So $\set{N_t}$ is simple symmetric at $z$, and lazy symmetric when not at $z$.
Let
$$ V'_t(z) = \sum_{j=0}^t \1{ N_j = z} , $$
be the number of times $\set{N_t}$ visits $z$ up to time $t$.

Define inductively $\rho_0 = \rho'_0 = 0$ and for $j \geq 0$,
$$ \rho_{j+1} = \min \set{ t \geq 1 \ : \ M_{\rho_{j}+t} = z } , $$
$$ \rho'_{j+1} = \min \set{ t \geq 1 \ : \ N_{\rho'_{j}+t} = z } . $$
If $N_t \geq M_t > z$ then
$$ \Pr \Br{ M_{t+1} = M_t + 1 }{ M_t } = \Pr \Br{ N_{t+1} = N_t + 1 }{ N_t } , $$
and
$$ \Pr \Br{ M_{t+1} = M_t - 1 }{ M_t } \geq \Pr \Br{ N_{t+1} = N_t - 1 }{ N_t } . $$
Thus, we can couple $M_{t+1}$ and $N_{t+1}$ so that $M_{t+1} \leq N_{t+1}$.
Similarly, if $N_t \leq M_t < z$ then $M_{t+1}$ moves towards $z$ with higher probability than
$N_{t+1}$, and they both move away from $z$ with probability $1/4$.
So we can couple $M_{t+1}$ and $N_{t+1}$ so that $M_{t+1} \geq N_{t+1}$.
If $N_t = M_t = z$ then $M_{t+1}$ and $N_{t+1}$ have the same distribution,
so they can be coupled so that $N_{t+1} = M_{t+1}$.

Thus, we can couple $\set{M_t}$ and $\set{N_t}$ so that for all $j \geq 0$,
$\rho_j \leq \rho_j'$ a.s.

Let $\set{S'_t}$ be a simple random walk on $\Z$.
For $x \in \Z$, let
$$ \tau_x = \min \set{ 2t \geq 2 \ : \ S'_{2t} = 2z \ , \ S'_0 = 2x } . $$
That is, $\tau_x$ is the first time a simple random walk started at $2x$ hits $2z$
(this is necessarily an even number).  In \cite[Chapter 9]{Revesz} it is shown that
$\tau_x$ has the same distribution as $\tau_{2z} - 2|z-x|$.
Note that if $N_t \neq z$ then
$S'_{2t+2} - S'_{2t}$ has the same distribution as $2(N_{t+1} - N_t)$.
Since $|N_{\rho'_{j-1} +1} -z|=1$, we get that for all $j \geq 2$,
$\rho'_j$ has the same distribution as $\frac{1}{2}(\tau_{2z}-2) +1$.
Also, $\rho'_1$ has the same distribution as $\frac{1}{2} \tau_0$ if $z \neq 0$, and
the same distribution as $\frac{1}{2} (\tau_{2z}-2) + 1$ if $z = 0$.
Hence, we conclude that for any $k \geq 1$,
$\sum_{j=1}^k \rho'_j$ has the same distribution as
$\frac{1}{2} \sum_{j=1}^k \tilde\rho_j$, where $\set{\tilde\rho_j}_{j \geq 1}$
are defined by
$$ \tilde\rho_{j+1} = \min \set{ 2t \geq 2 \ : \ S'_{\tilde\rho_j + 2t} = 2z } . $$
Finally note that $V_t(z) \geq k$ if and only if $\sum_{j=1}^k \rho_j \leq t$,
$V'_t(z) \geq k$ if and only if $\sum_{j=1}^k \rho'_j \leq t$,
and $L_t(2z) \geq k$ if and only if $\sum_{j=1}^k \tilde\rho_j \leq t$.
Thus, under the above coupling, for all $t \geq 0$, $V_t(z) \geq V'_t(z)$ a.s.
Also, $V'_t(z)$ has the same distribution as $L_{2t}(2z)$.
The lemma follows.
\end{proof}

\subsection{The Expectation of $X_t$}

Recall that $X_t = \sum_z X_t(z)$.

\begin{lem} \label{lem:E[Xt]}
There exists constants $c,C>0$ such that
for all $t \geq 0$,
$$ c t^{3/4} \leq \E [ X_t ] \leq C t^{3/4} . $$
\end{lem}

\begin{proof}
We first prove the upper bound.
For $z \in \Z$ let $A(z)$ be the indicator of the event
that the mixer reaches $z$ up to time $t$;  i.e. $A_t(z) = \1{V_t(z) \geq 1}$.
Note that $(\sigma_t(z)-z)(1-A_t(z)) = 0$.
Also, by definition $\sum_z V_t(z) = t$.
By Corollary \ref{cor:E[Xt(z)]},
using the Cauchy-Schwartz inequality,
\begin{align*}
\E [X_t] & = \sum_z \E[X_t(z)] \leq C_1 \sum_z \E [\sqrt{V_t(z)}] + 2 \E \sum_z A_t(z) \\
& \leq C_1 \E \sqrt{ \sum_z V_t(z) \cdot \sum_z A_t(z) } + 2 \E \sum_z A_t(z) ,
\end{align*}
for some constant $C_1>0$.
For any $z \in \Z$, if $A_t(z) = 1$, then there exists $0 \leq j \leq t$
such that $|S_j-z| = 1$.  That is, $A_t(z) = 1$ implies that $z \in [m_t-1,M_t+1]$,
where $M_t = \max_{0 \leq j \leq t} S_j$ and $m_t = \min_{0 \leq j \leq t} S_j$.
Thus, $\sum_z A(z) \leq M_t - m_t + 2$.  Since $M_t - m_t$ is just the number of sites
visited by a lazy random walk, we get (see e.g. \cite{Revesz})
$\E[\sum_z A_t(z)] \leq C_2 \sqrt{t}$, for some constant $C_2>0$.
Hence, there exists some constant $C_3>0$ such that
$$ \E[X_t] \leq C_1 \sqrt{t \cdot C_2 \sqrt{t}} + 2 C_2 \sqrt{t} \leq C_3 t^{3/4} . $$
This proves the upper bound.

We turn to the lower bound.
Let $\set{S'_t}$ be a simple random walk on $\Z$ started at $S'_0=0$,
and let
$$ L_t(z) = \sum_{j=0}^t \1{ S'_j=z } . $$
Let
$$ T(z) = \min \set{ t \geq 0 \ : \ S'_t = z } . $$
By the Markov property,
$$ \Pr \br{ L_{2t}(z) \geq k } \geq \Pr \br{ T(z) \leq t } \Pr \br{ L_t(0) \geq k } , $$
so
$$ \E [\sqrt{L_{2t}(2z)}] \geq \Pr \br{ T(2z) \leq t } \E [\sqrt{L_t(0)}] . $$
Theorem 9.3 of \cite{Revesz} can be used to show that
$\E[ \sqrt{L_t(0)} ] \geq c_1 t^{1/4}$,
for some constant $c_1>0$.
By Corollary \ref{cor:E[Xt(z)]},
and Lemma \ref{lem:distribution of Vt(z)},
there exists a constant $c_2>0$ such that
\begin{align*}
\E[X_t] & \geq c_2 \sum_z \E [\sqrt{L_{2t}(2z)}] - 2 \sum_z A_t(z) \\
& \geq c_1 t^{1/4} \cdot c_2 \E \sum_z \1{ T(2z) \leq t } - 2 C_2 \sqrt{t} .
\end{align*}
Let $M'_t = \max_{0 \leq j \leq t} S'_j$ and $m'_t = \min_{0 \leq j \leq t} S'_j$.
Then,
$$ \sum_z \1{T(2z) \leq t} = [m'_t,M'_t] \bigcap 2 \Z . $$
So for some constants $c_3,c_4>0$,
$$ \E [X_t] \geq c_3 t^{1/4} \cdot \frac{1}{2} \E[M'_t -m'_t-1] - 2 C_2 \sqrt{t}
\geq c_4 t^{3/4} . $$
\end{proof}

\section{Proof of Theorem \ref{thm:main thm}}
\label{sec:proof of thm}

\begin{proof}
Recall that $\Cov(z,\sigma)$ is the minimal length of a path on $\Z$,
started at $z$, that covers $\supp(\sigma)$.
Let $M_t = \max_{0 \leq j \leq t} S_j$ and $m_t = \min_{0 \leq j \leq t} S_j$,
and let $I_t = [m_t-1,M_t+1]$.
Note that $\supp(\sigma_t) \subset I_t$.  So
for any $z \in I_t$, $\Cov(z,\sigma_t) \leq 2|I_t|$.
$\set{S_t}$ has the distribution of a lazy random walk on $\Z$,
so $\set{2 S_t}$ has the same distribution as $\set{S'_{2t}}$, where
$\set{S'_t}$ is a simple random walk on $\Z$.  It is well known (see e.g. \cite[Chapter 2]{Revesz})
that there exist constants $c_1,C_1>0$ such that $c_1 \sqrt{t} \leq \E[|I_t|] \leq C_1 \sqrt{t}$.
Since $S_t \in I_t$, we get that $\E [\Cov(S_t,\sigma_t)] \leq 2 C_1 \sqrt{t}$.
Together with Propositions \ref{prop:distance upper in mixer group} and
\ref{prop:distance lower in mixer group}, and with Lemma \ref{lem:E[Xt]}, we get that
there exist constants $c,C>0$ such that for all $t \geq 0$,
$$ c t^{3/4} \leq \frac{1}{2} \E[X_t] \leq \E[D_t] \leq 2 \E[\Cov(S_t,\sigma_t)] + 5 \E[X_t] \leq C t^{3/4} . $$
\end{proof}

%
%
%
%


\bibliographystyle{acm}

\end{document}